\documentclass[twocolumn,showpacs,preprintnumbers,amsmath,amssymb]{revtex4}

\usepackage{graphicx}
\usepackage{dcolumn}
\usepackage{bm}
\usepackage{epsfig}
\usepackage{subfigure}

\begin{document}

\title{Optimal Design of Minimum Energy Pulses for Bloch Equations\\ in the case of Dominant Transverse Relaxation}

\author{Dionisis Stefanatos}
\email{dionisis@post.harvard.edu}
\affiliation{Prefecture of Kefalonia, Argostoli, Kefalonia 28100, Greece}

\date{\today}

\begin{abstract}
In this report, we apply Optimal Control Theory to design minimum energy $\pi/2$ and $\pi$ pulses for Bloch equations, in the case where transverse relaxation rate is much larger than longitudinal so the later can be neglected. Using Pontryagin's Maximum Principle, we derive an optimal feedback law and subsequently use it to obtain analytical expressions for the energy and duration of the optimal pulses.
\end{abstract}

\maketitle

Optimal Control Theory \cite{Pontryagin} has been extensively used recently for the design of pulses that optimize the performance of various Nuclear Magnetic Resonance (NMR) and Quantum Optics systems limited by the presence of relaxation \cite{Khaneja03_1,Khaneja03_2,BBCROP,TROPIC,Stefanatos04,Stefanatos05,Sklarz,Sugny,Li}. In this report we use it to derive minimum energy $\pi/2$ and $\pi$ pulses for Bloch equations, in the case where transverse relaxation dominates.

\begin{figure}[h]
\centering
\includegraphics[scale=0.35]{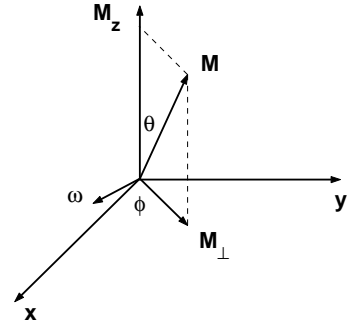}
\caption{The optimal transverse magnetic field $\omega$ is perpendicular to $M_\perp$ and its phase is constant. Without loss of generality, the experimental setup can be arranged such that $\omega \parallel x$-axis. In this case, $\phi=\pi/2$ and $\mathbf{M}$ rotates in $yz$-plane.}
\label{fig:spherical}
\end{figure}

The Bloch equations, in a resonant rotating frame and when longitudinal relaxation is neglected are \cite{Ernst}
\begin{eqnarray}
	\dot{M}_z & = & \omega_yM_x-\omega_xM_y \nonumber\\
	\dot{M}_x & = & -RM_x-\omega_yM_z \nonumber\\
    \dot{M}_y & = & -RM_y+\omega_xM_z \nonumber
\end{eqnarray}
where $\mathbf{M}=(M_x, M_y, M_z)$ is the magnetization vector, $\omega_x,\omega_y$ are the transverse components of the magnetic field and $R>0$ is the transverse relaxation rate. If we make the following change of variables, see Fig. \ref{fig:spherical}
\begin{eqnarray}
	a & = & \ln[M(t)/M(0)] \nonumber\\
	\tan\theta & = & \sqrt{M^2_x+M^2_y}/M_z \nonumber\\
    \tan\phi & = & M_y/M_x \nonumber
\end{eqnarray}
where $M = \sqrt{M^2_x+M^2_y+M^2_z}$, we obtain
\begin{eqnarray}
	\label{magnitude}\dot{a} & = & -R\sin^2\theta \\
	\label{angle}\dot{\theta} & = & \omega_\perp-R\sin\theta\cos\theta \\
    \label{phi}\dot{\phi} & = & \omega_\parallel\cot\theta
\end{eqnarray}
where $\omega_\perp = \omega_x\sin\phi-\omega_y\cos\phi, \omega_\parallel = \omega_x\cos\phi+\omega_y\sin\phi$ are the components of transverse magnetic field perpendicular and parallel to $\mathbf{M}_\perp = (M_x,M_y)$, respectively. Note that $\omega_\parallel$ does not affect the angle $\theta$ of the pulse. It just rotates $\mathbf{M}$ around $z$-axis, resulting in a waste of energy. Thus, optimality requires $\omega_\parallel=0\Rightarrow \phi=\mbox{constant}$. Equations (\ref{magnitude}) and (\ref{angle}) are sufficient to describe the rotation and from now on we use $\omega$ to denote $\omega_\perp$, see Fig. \ref{fig:spherical}.

Observe that when $\omega$ is unbounded, we can rotate $\theta$ instantaneously to the desired final value $\theta(T)=\pi/2\,\mbox{or}\,\pi$ without losses in $a$, i.e. with $a(T)=a(0)=0$ (equivalently $M(T)=M(0)$). This transfer requires an infinite amount of energy so it is unrealistic. A meaningful optimization problem is the following: For a specified final value $a(T)<a(0)=0$ (equivalently specified $M(T)<M(0)$), what is the optimal control $\omega(t)$, $0\leq t\leq T$, that accomplishes the transfer $(a(0)=0,\theta(0)=0)\rightarrow(a(T),\theta(T)=\pi/2\,\mbox{or}\,\pi)$, while minimizing energy $\int_0^T\omega^2(t)/2dt$? The control Hamiltonian \cite{Pontryagin} for this problem is
\begin{equation}
\label{hamiltonian}
H = -\omega^2/2+\lambda_{\theta}(\omega-R\sin\theta\cos\theta)-\lambda_a R\sin^2\theta
\end{equation}
where $\lambda_{\theta},\lambda_a$ are the Lagrange multipliers. According to Pontryagin's maximum principle \cite{Pontryagin}, necessary conditions for optimality of $(\omega(t), a(t), \theta(t), \lambda_{\theta}(t),\lambda_a(t))$ are
\begin{eqnarray}
\label{optimality}
\vartheta H/\vartheta\omega & = & 0 \Rightarrow\omega = \lambda_{\theta}\\
\label{adjoint1}
\dot{\lambda}_\theta & = & -\vartheta H/\vartheta\theta = \lambda_\theta R\cos 2\theta+\lambda_a R\sin 2\theta\\
\label{adjoint2}
\dot{\lambda}_a & = & -\vartheta H/\vartheta a = 0 \Rightarrow \lambda_a=\mbox{constant}
\end{eqnarray}
Additionally, the optimal $(\omega, a, \theta, \lambda_{\theta},\lambda_a)$ satisfies \cite{Pontryagin}
\begin{equation}
H(\omega,a,\theta,\lambda_\theta,\lambda_a)=0,\; 0\leq t\leq T.
\end{equation}
Using (\ref{hamiltonian}),(\ref{optimality}) the above condition becomes
\begin{equation}
\label{optimallambda}
\lambda^2_\theta-2\lambda_\theta R\sin\theta\cos\theta-2\lambda_a R\sin^2\theta=0
\end{equation}
Note that when $\theta(t_{\pi/2})=\pi/2$ (there is such a time $t_{\pi/2}$ for both final values $\theta(T)$ that we consider here), then (\ref{optimallambda}) gives $\lambda_a=\lambda^2_\theta(t_{\pi/2})/2R\geq 0$. Solving the above quadratic equation for $\lambda_\theta$ and using (\ref{optimality}), we find the optimal control $\omega$. Note that only the positive solution of the quadratic equation has physical meaning (corresponds to increasing $\theta$) for the cases that we study here. The optimal $\omega$ is given by the following feedback law
\begin{equation}
\label{omega}
\omega(\theta)=R\sin\theta(\cos\theta+\sqrt{\cos^2\theta+\kappa^2})
\end{equation}
where $\kappa^2=2\lambda_a/R=\lambda^2_\theta(t_{\pi/2})/R^2$. Using (\ref{omega}), the validity of (\ref{adjoint1}) can be easily verified.
Inserting (\ref{omega}) in (\ref{angle}) we obtain the differential equation for the optimal trajectory
\begin{equation}
\label{theta}
\dot{\theta}=R\sin\theta\sqrt{\cos^2\theta+\kappa^2}
\end{equation}
Integrating (\ref{magnitude}) from $t=0$ to the final time $t=T$ we get
\begin{equation}
\label{constraint}
\int_0^T\sin^2\theta dt = -\frac{a(T)}{R} = \frac{1}{R}\ln\frac{M(0)}{M(T)}
\end{equation}
If we use (\ref{theta}) to change the integration from time to angle, we obtain
\begin{equation}
\int_{\theta(0)}^{\theta(T)}\frac{\sin\theta}{\sqrt{\cos^2\theta+\kappa^2}}\,d\theta=-\ln r
\end{equation}
where $r=M(T)/M(0)<1$. This condition determines $\kappa$ and the results for the $\pi/2$ and $\pi$ pulses are
\begin{equation}
\kappa_{\pi/2}=\frac{2r}{1-r^2},\;\;\kappa_{\pi}=\frac{2\sqrt{r}}{1-r}
\end{equation}

The duration of the optimal pulses is determined from the relation
\begin{equation}
T=\int_0^Tdt=\int_{\theta(0)}^{\theta(T)}\frac{1}{\dot{\theta}(\theta)}\,d\theta
\end{equation}
using again (\ref{theta}). Note that since $\theta=0,\pi$ are equilibrium points for (\ref{theta}), we actually start from a small positive initial value $\theta(0)=\epsilon$ for both cases, and additionally for the $\pi$ pulse we end to the value $\theta(T)=\pi-\epsilon$. The duration of the optimal pulses as $\epsilon\rightarrow 0$ is
\begin{eqnarray}
T_{\pi/2} & = & \frac{1}{R}\frac{1-r^2}{1+r^2}\left[\ln\left(\frac{1+r^2}{r}\right)-\ln \epsilon\right]\\
T_{\pi} & = & \frac{1}{R}\frac{1-r}{1+r}\left[\ln\left(\frac{(1+r)^2}{r}\right)-2\ln \epsilon\right]
\end{eqnarray}

Finally, the energy of the optimal pulses is calculated from
\begin{equation}
\int_0^T\frac{\omega^2(t)}{2}dt=\int_{\theta(0)}^{\theta(T)}\frac{\omega^2(\theta)}{2\dot{\theta}(\theta)}\,d\theta
\end{equation}
using (\ref{omega}) and (\ref{theta}). The result is
\begin{equation}
E_{\pi/2}=R\frac{1}{1-r^2},\;\;E_{\pi}=R\frac{1+r}{1-r}
\end{equation}
Observe that for $r\rightarrow 1$ we have $T_{\pi/2}, T_{\pi}\rightarrow 0$ and $E_{\pi/2}, E_{\pi}\rightarrow \infty$, as mentioned above.

\begin{figure}[t]
 \centering
		\begin{tabular}{cc}
     	\subfigure[$\ $Optimal $\pi/2$ pulse]{
	            \label{fig:pi2}
	            \includegraphics[width=.45\linewidth]{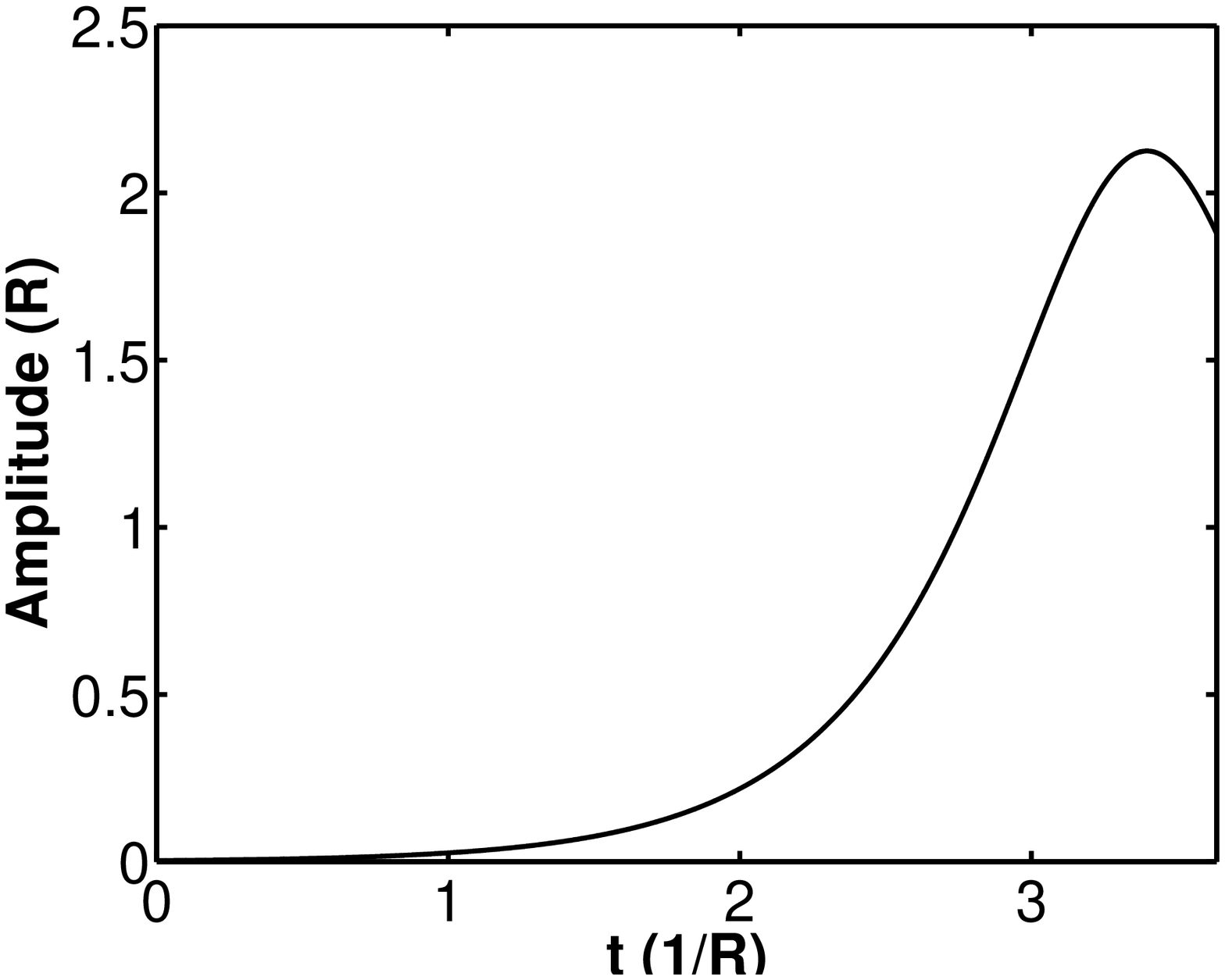}} &
	        \subfigure[$\ $Optimal $\pi$ pulse]{
	            \label{fig:pi}
	            \includegraphics[width=.45\linewidth]{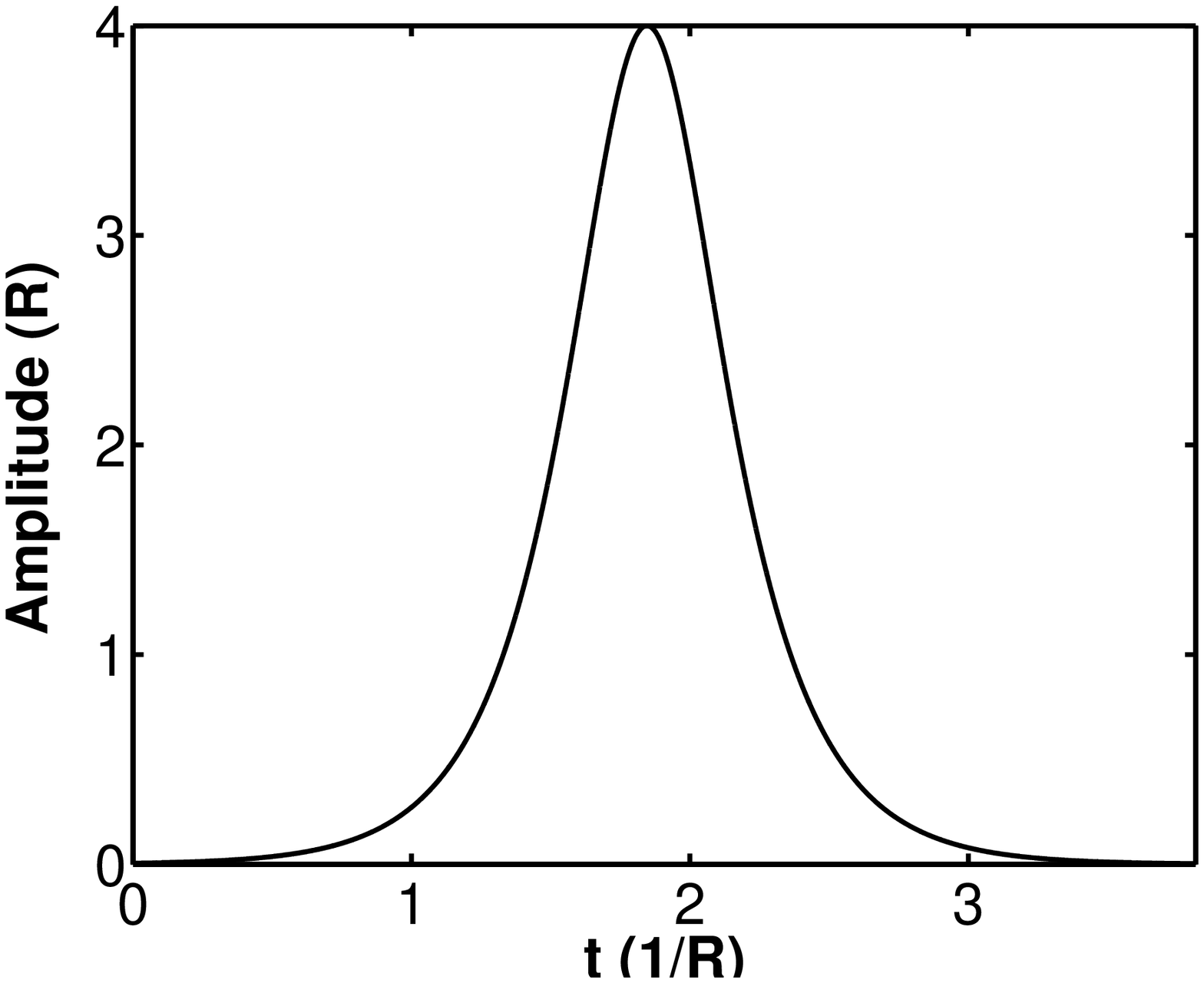}} \\
	        \subfigure[$\ $Optimal trajectory]{
	            \label{fig:pi2_traj}
	            \includegraphics[width=.45\linewidth]{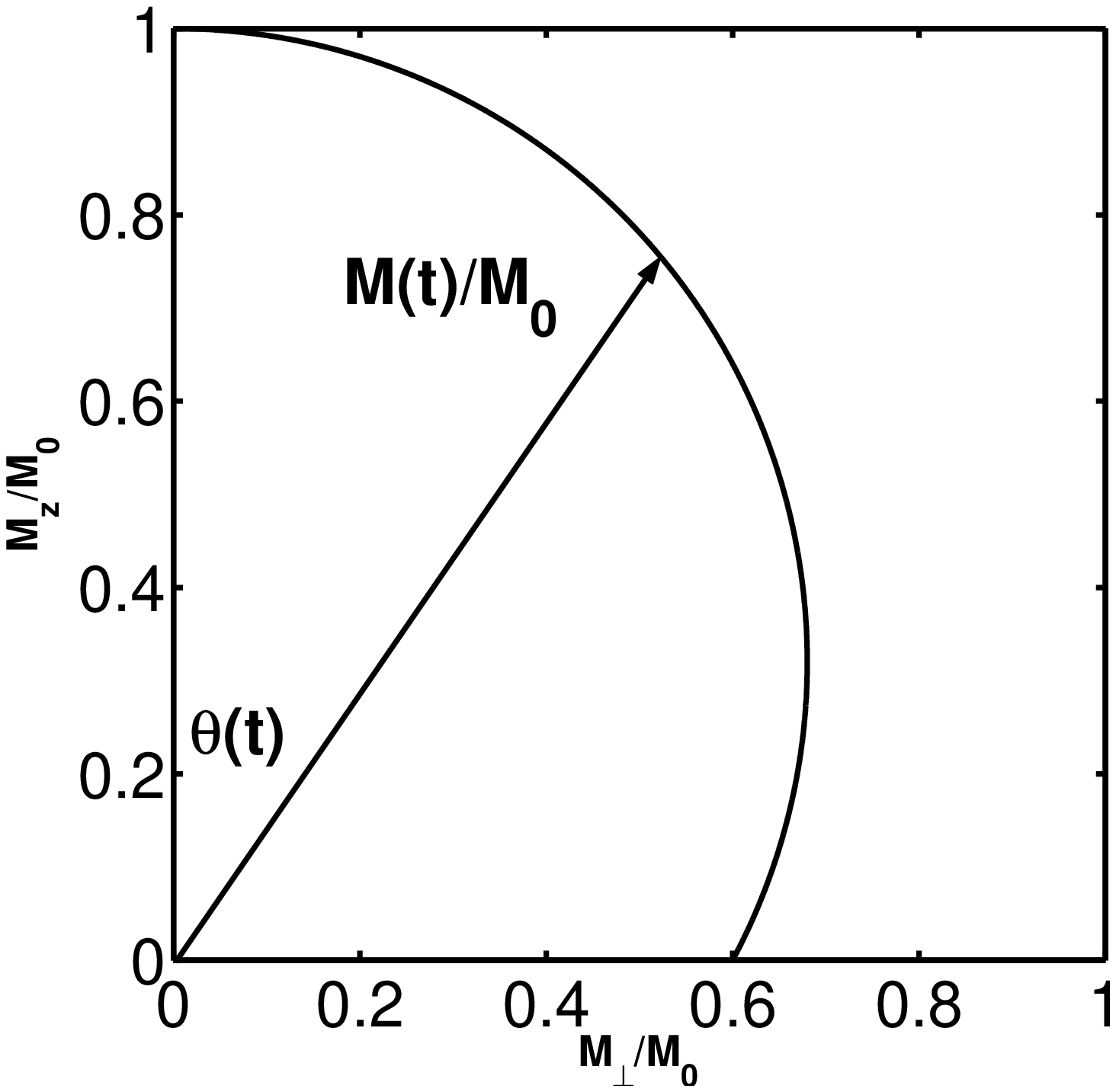}} &
			\subfigure[$\ $Optimal trajectory]{
	            \label{fig:pi_traj}
	            \includegraphics[width=.45\linewidth]{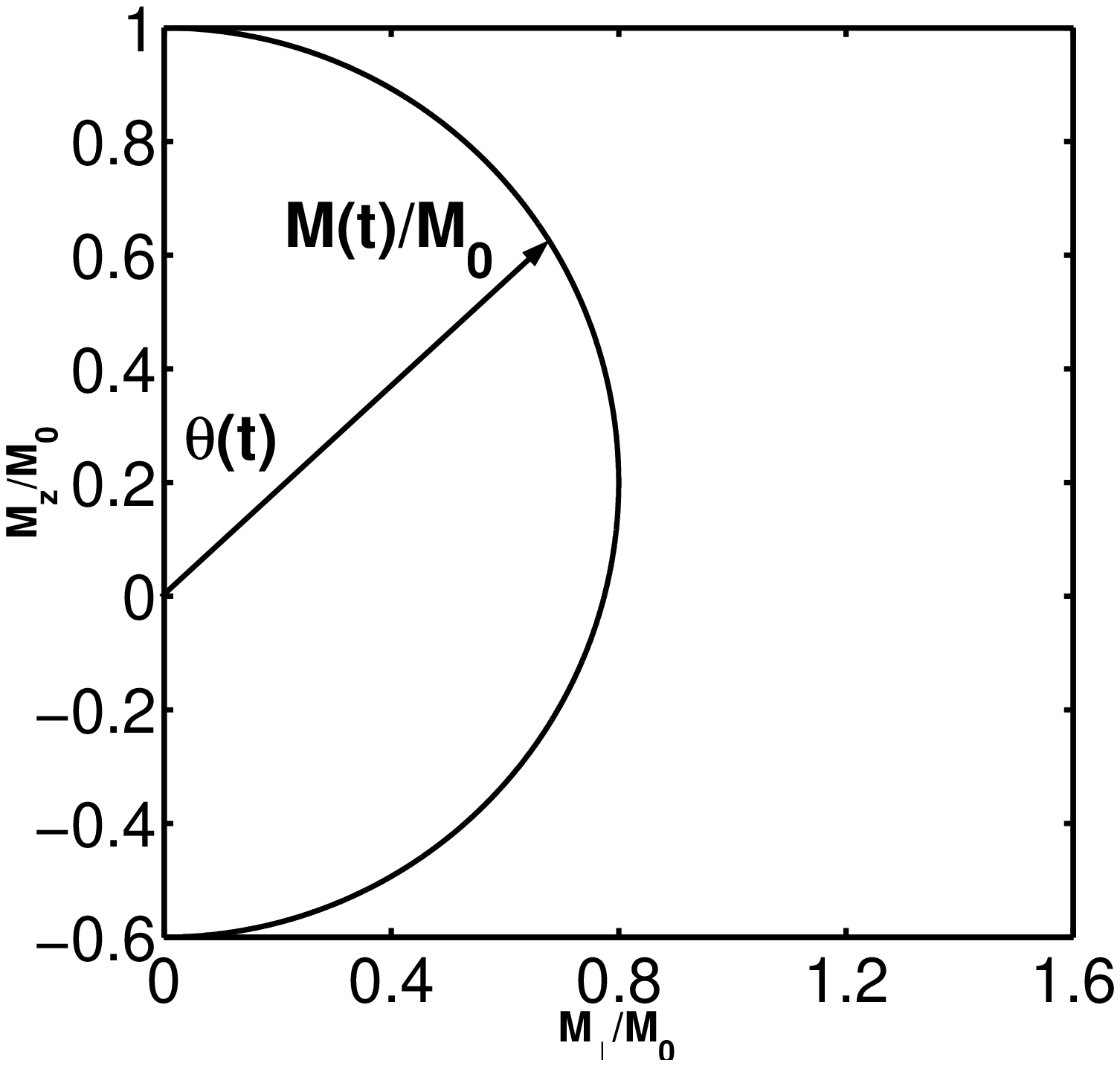}} \\
		\end{tabular}
 \caption{Minimum energy $\pi/2$ (panel a) and $\pi$ (panel b) pulses for $M(0)=M_0, M(T)=0.6M_0, \epsilon=10^{-3}$. The corresponding trajectories are also shown (panels c,d).}
 \label{fig:pulses}
\end{figure}

\begin{figure}[t]
\centering
\includegraphics[scale=0.35]{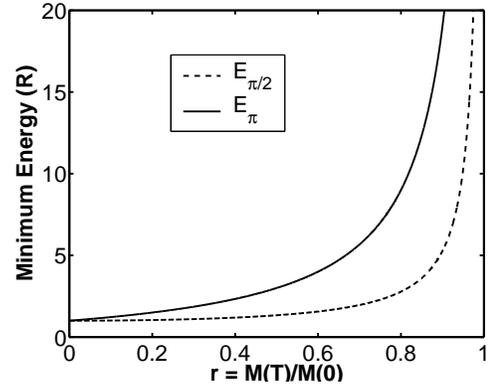}
\caption{The energy of the optimal pulses as a function of the ratio $r=M(T)/M(0)$.}
\label{fig:energyplot}
\end{figure}

In Fig. \ref{fig:pulses} we plot the optimal $\pi/2$ and $\pi$ pulses for $M(0)=M_0, M(T)=0.6M_0, \epsilon=10^{-3}$, as well as the corresponding trajectories of normalized magnetization vector. In Fig. \ref{fig:energyplot} we plot the energy of the pulses as a function of the ratio $r=M(T)/M(0)$.


To conclude, in this report we calculated minimum energy $\pi/2$ and $\pi$ pulses for Bloch equations in the case where transverse relaxation dominates, using Optimal Control Theory. We expect this analytical work to serve as a reference for numerical studies of more complicated and realistic situations that incorporate for example longitudinal relaxation and magnetic field inhomogeneity.

\noindent\textbf{Acknowledgements}: The author would like to thank Jr-Shin Li for his constant support.


\begin{thebibliography}{99}


\bibitem{Pontryagin}
L.S. Pontryagin, V.G. Boltyanskii, R.V. Gamkrelidze, and E.F.
  Mishchenko, {\it The Mathematical Theory of Optimal Processes} (Interscience Publishers, New York, 1962).

\bibitem{Khaneja03_1}
N. Khaneja, T. Reiss, B. Luy, S. J. Glasser, J. Magn. Reson. 162 311 (2003).

\bibitem{Khaneja03_2}
N. Khaneja, B. Luy, and S. J. Glaser, Proc. Natl. Acad. Sci. USA 100, 13162 (2003).

\bibitem{BBCROP}
N. Khaneja, J.-S. Li, C. Kehlet, B. Luy and S. J. Glaser, Proc. Natl. Acad. Sci. USA 101, 14742 (2004).

\bibitem{TROPIC}
D. P. Frueh, T. Ito, J.-S. Li, G. Wagner, S. J. Glaser and N. Khaneja, J. of Biomol. NMR 32, 23 (2005).

\bibitem{Stefanatos04}
D. Stefanatos, N. Khaneja, and S. J. Glaser, Phys. Rev. A 69, 022319 (2004).

\bibitem{Stefanatos05}
D. Stefanatos, S. J. Glaser, and N. Khaneja, Phys. Rev. A 72, 062320 (2005).

\bibitem{Sklarz}
S.E. Sklarz, D.J. Tannor, and N. Khaneja, Phys. Rev. A 69, 053408 (2004).

\bibitem{Sugny}
D. Sugny, C. Kontz, and H. R. Jauslin, Phys. Rev. A 76, 023419 (2007).

\bibitem{Li}
J.-S. Li, J. Ruths, and D. Stefanatos, arXiv:0908.2093v1 [physics.chem-ph].

\bibitem{Ernst}
R. R. Ernst, G. Bodenhausen, A. Wokaun, {\it Principles of Nuclear Magnetic Resonance in One and Two Dimensions} (Clarendon Press, Oxford, 1987).

\end{thebibliography}

\end{document}